\documentclass[12pt,reqno,a4paper]{amsart} 
\usepackage[headings]{fullpage}
\usepackage{hyperref} 
\usepackage{datetime}
\hypersetup{breaklinks=true,  colorlinks=false}   

\usepackage[english]{babel} 
 
\renewcommand{\qedsymbol}{$\square$}
\newenvironment{Proof}[1][Proof]{\par\noindent\textbf{#1.}~}
{\hfill\qedsymbol\smallskip\par}

\newcommand{\dx}{\mathrm{d}}
\newcommand{\eps}{\varepsilon}
\newcommand{\R}{\mathbb{R}}
\newcommand{\N}{\mathbb{N}} 

\newcommand{\C}{\mathbb{C}}

\newcommand{\Stilde}{\widetilde{S}}
\newcommand{\Etilde}{\widetilde{E}}

\newcommand{\second}{\prime\prime}

\newcommand{\Odip}[2]{\mathcal{O}_{#1}\!\left(#2\right)\mathchoice{\!}{}{}{}}
\newcommand{\Odig}[1]{\mathcal{O}\Bigl(#1\Bigr)\mathchoice{\!}{}{}{}}

\newcommand{\Odi}[1]{\Odip{}{#1}}
\newcommand{\odip}[2]{{o}_{#1}\!\left(#2\right)\mathchoice{\!}{}{}{}}
\newcommand{\odi}[1]{\odip{}{#1}}

\newtheoremstyle{sltheorems}
{10pt}
{6pt}
{\slshape}
{}
{\bfseries}
{.}
{.5em}
{\thmname{#1}\thmnumber{ #2}\thmnote{ (#3)}}

\theoremstyle{sltheorems}
\newtheorem{Theorem}{Theorem}
\newtheorem{Lemma}{Lemma}

\allowdisplaybreaks

\begin{document}

\title[one prime and two squares of primes in short intervals]{Sum of one prime and two squares of primes \\ in short intervals} 
\author[]{Alessandro Languasco \lowercase{and} Alessandro Zaccagnini}

\date{\today, \currenttime}
\begin{abstract}
Assuming the Riemann Hypothesis we prove that the interval $[N, N + H]$
contains an integer which is a sum of a prime and two squares of
primes provided that $H \ge C (\log N)^{4}$, where $C > 0$ is an
effective constant.
\end{abstract}

\selectlanguage{english}
\maketitle

\section{Introduction}
The problem of representing an integer as a sum of a prime and of two prime squares
is classical. 
Letting
\[
{\mathcal A} = \{ n\in \N : n\equiv 1 \bmod 2;\ n \not \equiv 2\bmod 3\},
\]
it is conjectured that every sufficiently large $n\in {\mathcal A}$  can be represented as
$n= p_{1} +p_{2}^{2} + p_{3}^{2}$. Let now $N$ be a large integer.
Several results about the cardinality $E(N)$ of the set
of integers $n\leq N$, $n \in {\mathcal A}$ which are not representable as a sum of a prime and two prime squares
were proved during the last 75 years; we recall the papers of  Hua \cite{Hua1938}, Schwarz \cite{Schwarz1961}, Leung-Liu \cite{LeungL1993}, Wang \cite{Wang2004}, Wang-Meng \cite{WangM2006a}, Li \cite{Li2008} and Harman-Kumchev \cite{HarmanK2010}.
Recently L.~Zhao \cite{Zhao2014a} proved that 
\begin{equation*}
E(N)\ll N^{1/3+\eps}.
\end{equation*} 
As a consequence we can say that every integer  $n\in [1,N]\cap  {\mathcal A}$, with at most $\Odi{N^{1/3+\eps}}$ exceptions,
 is the sum of a prime and two prime squares. Letting
\begin{equation}
\label{r-def}
  r(n)
  =
  \sum_{p_{1} + p_{2}^2 + p_{3}^{2} = n} \log p_{1} \log p_{2} \log p_{3},
\end{equation} 
in fact L.~Zhao also proved that a suitable asymptotic formula
for $r(n)$ holds for every $n\in [1,N]\cap  {\mathcal A}$, with at most $\Odi{N^{1/3+\eps}}$ exceptions.

In this paper we study the average behaviour  of $r(n)$  
over  short intervals $[N,N+H]$, $H=\odi{N}$.
Assuming that the Riemann Hypothesis (RH) holds,
we prove that a suitable asymptotic formula for such an average of
$r(n)$ holds in short intervals with no exceptions.

\begin{Theorem}
\label{existence}
Assume the Riemann Hypothesis (RH). We have
 \[
   \sum_{n = N+1}^{N + H}   r(n) = \frac{\pi}{4} HN
   +
\Odig{H^{1/2}N (\log N)^{2} +H N^{3/4}(\log N)^{3}   +  H^{2}L^{3/2}}
    \quad \textrm{as} \ N \to \infty,
 \] 
uniformly for  $\infty((\log N)^{4})\le H \le \odi{NL^{-3/2}}$, where $f=\infty(g)$ means $g=\odi{f}$. 
\end{Theorem}

Letting
 \[
  r^{*}(n)
  =
  \sum_{p_{1} + p_{2}^2 + p_{3}^{2} = n} 1,
\]
a similar asymptotic formula holds for the average of $r^{*}(n)$ too.

In the unconditional case our proof yields a weaker result
than Zhao's, namely, the asymptotic formula for the average of $r(n)$
holds just for $H\ge N^{7/12+\eps}$; 
for this reason, here we are only concerned with the conditional one.
It is worth remarking that, under the assumption of RH, 
the formula in Theorem \ref{existence} implies that
every interval $[N,N+H]$ contains an integer which is a sum of a prime 
and two prime squares, where $C L^{4} \le  H =\odi{NL^{-3/2}}$,
$C>0$ is a suitable large constant and $L = \log N$.
We recall that the  analogue results for the binary Goldbach problem are respectively 
$H\gg N^{c+\eps}$ with $c=21/800$, by Baker-Harman-Pintz  and Jia, see \cite{Pintz2012},
and $H\gg L^{2}$, under the assumption of RH; see, \emph{e.g.}, \cite{LanguascoP1994}. 
Assuming RH, the expectation in Theorem \ref{existence} is the lower
bound $H \gg L^{2}$ since the crucial error term should be
$\ll H^{1/2}N \log N$; the loss of
a factor $L$ in such an error term is due to the lack of information about a truncated fourth-power average for $\Stilde_{2}(\alpha)$: see Lemma \ref{Hua-Rieger-Lemma} and \eqref{I3-estim-final} below.

The proof of Theorem \ref{existence} uses the original
Hardy-Littlewood settings of the circle method,
\emph{i.e.}, with infinite series instead of finite sums over primes. 
This is due to the fact that for this problem both the direct and the finite sums 
approaches do not seem to be able to work in intervals shorter than $N^{1/2}$.

\medskip
{\bf Acknowledgements.}    This research was partially supported by the grant PRIN2010-11 \textsl{Arithmetic Algebraic Geometry and Number Theory}. 
We wish to thank the referee for pointing out some inaccuracies.

\section{Notation and Lemmas}

Let $\ell \ge 1$ be an integer. The standard circle method
approach requires to define
\begin{equation*}
S_{\ell}(\alpha)
 =
\sum_{1 \le p^{\ell} \le  N} \log p \ e(p^{\ell} \alpha )
\quad
\textrm{and}
\quad
  T_{\ell}(\alpha)
  =
    \sum_{1 \le n^{\ell} \le  N} e(n^{\ell} \alpha ),
\end{equation*}
where $e(x)=\exp(2\pi i x)$,  and
needs the following lemma which collects the results 
of Theorems 3.1-3.2 of \cite{LanguascoZ2013b}.
\begin{Lemma}  
\label{App-BCP-Gallagher}
Let $N$ be a large integer, $\ell  > 0$ be a real number and $\eps$ be an arbitrarily small
positive constant. Then there exists a positive constant 
$c_1 = c_{1}(\eps)$, which does not depend on $\ell$, such that
\[
\int_{-1/H}^{1/H}
\vert
S_{\ell}(\alpha) - T_{\ell}(\alpha)  
\vert^2 
\, \dx \alpha
\ll_{\ell}   N^{2/\ell -1}
\Bigl(
\exp \Big( - c_{1}  \Big( \frac{L}{\log L} \Big)^{1/3} \Big)
+
\frac{H L^{2}}{N}
\Bigr),
\]
uniformly for  $N^{1-5/(6\ell)+\eps}\le H \le N$.  
Assuming further RH we get 
\[
\int_{-1/H}^{1/H}
\vert
S_{\ell}(\alpha) - T_{\ell}(\alpha) 
\vert^2 
\, \dx \alpha
\ll_{\ell} 
\frac{N^{1/\ell} L^{2}}{H} + H N^{2/\ell-2} L^{2},
\]
uniformly for  $N^{1-1/\ell}\le H \le N$.  
\end{Lemma}
So it is clear that this approach works only when  the lower bound  $H\ge N^{1-1/\ell}$
holds.
Such a limitation comes from the fact that Gallagher's lemma translates the 
mean-square average of an exponential sum in a short interval problem. 
When  $\ell$-powers are involved, this leads to $p^{\ell}\in[N,N+H]$ 
which is a non-trivial condition only when $H\ge N^{1-1/\ell}$. 

So, when $\ell=2$, the standard circle method approach works
only if $H\ge N^{1/2}$; on the other hand we can easily show that the
direct attack works, under RH, only for $H=\infty(N^{1/2}L^2)$.
Therefore, to have the chance to reach smaller $H$-values, 
we will use the original Hardy and 
Littlewood \cite{HardyL1923} circle  method setting, \emph{i.e.},
the weighted exponential sum 
\begin{equation*}
\Stilde_\ell(\alpha)
=
\sum_{n=1}^{\infty} 
\Lambda(n) e^{-n^{\ell}/N}
e(n^{\ell}\alpha),
\end{equation*}
since it lets us avoid the use of Gallagher's lemma, see 
Lemmas \ref{Linnik-lemma}-\ref{LP-Lemma-gen} below.

The first ingredient we need is the following explicit formula
which generalizes and slightly sharpens what Linnik \cite{Linnik1946} proved:
see also eq.~(4.1) of \cite{Linnik1952}.
Let 
\begin{equation}
\label{z-def}
z= 1/N-2\pi i\alpha.
\end{equation}
We remark that 

\begin{equation}
\label{z-estim}
\vert z\vert ^{-1} \ll \min \bigl(N, \vert \alpha \vert^{-1}\bigr).
\end{equation}

\begin{Lemma} 
\label{Linnik-lemma}
Let $\ell \ge 1$ be an integer, $N \ge 2$  and $\alpha\in [-1/2,1/2]$.
Then
\begin{equation}
\label{expl-form}
\Stilde_{\ell}(\alpha)  
= 
\frac{\Gamma(1/\ell)}{\ell z^{1/\ell}}
- 
\frac{1}{\ell}\sum_{\rho}z^{-\rho/\ell}\Gamma\Bigl(\frac{\rho}{\ell}\Bigr) 
+
\Odip{\ell}{1},
\end{equation}
where $\rho=\beta+i\gamma$ runs over
the non-trivial zeros of $\zeta(s)$.
\end{Lemma} 

\begin{Proof}
We recall that Linnik proved this formula in the case $\ell =1$, with an error term $\ll 1+\log^3 (N \vert\alpha\vert)$.

Following the line of Lemma 4 in Hardy and Littlewood \cite{HardyL1923}
and  of \S4 in Linnik \cite{Linnik1946},  we have that 
\begin{equation}
\label{Mellin}
\Stilde_\ell(\alpha)  
= 
\frac{\Gamma(1/\ell)}{\ell z^{1/\ell}}
- 
\frac{1}{\ell}\sum_{\rho}z^{-\rho/\ell}\Gamma\Bigl(\frac{\rho}{\ell}\Bigr) 
-
\frac{\zeta'}{\zeta}(0)
-
\frac{1}{2\pi i}
\int_{(-\sqrt{3}/2)} 
\frac{\zeta'}{\zeta}(\ell w) \Gamma(w)z^{-w} \, \dx w.
\end{equation}
Now we estimate the integral in \eqref{Mellin}.  
Writing $w=-\sqrt{3}/2+it$, we have
$\vert (\zeta'/\zeta)(\ell w)\vert \ll_{\ell} \log (\vert t \vert +2)$,
$z^{-w}= \vert z \vert^{\sqrt{3}/2} \exp(t \arg(z))$,
where $\vert \arg(z) \vert \le  \pi/2$.
Furthermore the Stirling formula implies that
$\Gamma(w) \ll \vert t \vert^{-(\sqrt{3}+1)/2} \exp(-\pi\vert t \vert/2)$.
Hence
\begin{align*}
  \int_{(-\sqrt{3}/2)}
    \frac{\zeta'}{\zeta}(\ell w) \Gamma(w) z^{-w} \, \dx w
  &\ll_{\ell}
    \vert z \vert^{\sqrt{3}/2}
  \int_0^1 \log(t + 2) \, \dx t \\
  &\qquad
  +
    \vert z \vert^{\sqrt{3}/2}
  \int_1^{\infty}
    \log (t + 2) t^{-(\sqrt{3}+1)/2}
    \exp\Bigl( (\arg(z) - \frac{\pi}{2}) t \Bigr) \, \dx t \\
  &\ll_{\ell}
    \vert z \vert^{\sqrt{3}/2}  
  +
    \vert z \vert^{\sqrt{3}/2}
  \int_1^{\infty}
    \log(t + 2) t^{-(\sqrt{3}+1)/2} \, \dx t 
\ll_{\ell}
  \vert z \vert^{\sqrt{3}/2}.
\end{align*}
This is $\ll_{\ell} 1 $ as stated since  $z \ll 1$ by \eqref{z-def}.
Hence the lemma is proved.
\end{Proof}
We explicitly remark that Lemma~\ref{Linnik-lemma} is stronger than the corresponding Lemma 1 of \cite{LanguascoZ2014a} (or Lemma 1 of \cite{LanguascoZ2013a}) because in this case $\alpha$ is bounded.

The second lemma is an $L^{2}$-estimate of the 
remainder term in \eqref{expl-form} 
which generalizes a result of  Languasco and Perelli \cite{LanguascoP1994};
we will follow their proof inserting many details since the presence of $\ell$ 
changes the shape of the involved estimates at several places.
In fact we will use Lemma \ref{LP-Lemma-gen} just for $\ell=1,2$ but we 
take this occasion to describe the more general case  since it may be 
useful for future works.
\begin{Lemma}
\label{LP-Lemma-gen}
Assume RH. Let $\ell \ge 1$ be an integer and $N$ be a
sufficiently large integer. For  $0 \le \xi \le 1/2$, we have 
\[
\int_{-\xi}^{\xi} \,
\Bigl\vert
\Stilde_\ell(\alpha) - \frac{\Gamma(1/\ell)}{\ell z^{1/\ell}}
\Bigr\vert^{2}
\dx \alpha 
\ll_{\ell}
N^{1/\ell}\xi L^{2}.
\]
\end{Lemma} 

\begin{Proof}
Since 
$z^{-\rho/\ell} =  \vert z \vert^{-\rho/\ell} \exp\bigl(-i(\rho/\ell)\arctan2\pi N\alpha\bigr)$,
by RH and Stirling's formula we have that
\[
\frac{1}{\ell}
\sum_{\rho}z^{-\rho/\ell}\Gamma\Bigl(\frac{\rho}{\ell}\Bigr) 
\ll_{\ell}
\sum_{\rho}  \vert z \vert^{-1/(2\ell)}
\vert \gamma \vert ^{(1-\ell)/(2\ell)}
\exp
\Bigl(\frac{\gamma}{\ell}\arctan2\pi N\alpha - \frac{\pi}{2\ell} \vert\gamma \vert\Bigr).
\]
If $\gamma\alpha \le 0$ or $ \vert\alpha \vert\leq\ 1/N$ we get
\(
\sum_{\rho}z^{-\rho/\ell}\Gamma (\rho/\ell) 
\ll_{\ell} N^{1/(2\ell)} ,
\)
where, in the first case, $\rho$ runs over the zeros with $\gamma\alpha \le 0$.
Hence 
\begin{equation}
\label{R1-first}
I(N,\xi,\ell)
:=
\int_{-\xi}^{\xi} \,
\Bigl\vert
\Stilde_\ell(\alpha) - \frac{\Gamma(1/\ell)}{\ell z^{1/\ell}}
\Bigr\vert^{2} \dx \alpha 
\ll_{\ell} 
N^{1/\ell}\xi
\end{equation}
if $0\leq\xi\le 1/N$, and
\begin{equation}
\label{R1-estim}
I(N,\xi,\ell)
\ll_{\ell}
\int_{1/N}^{\xi}
\Big \vert\sum_{\gamma>0}z^{-\rho/\ell}\Gamma\Bigl(\frac{\rho}{\ell}\Bigr) \Big \vert^2 \dx \alpha +
\int_{-\xi}^{-1/N}
\Big \vert\sum_{\gamma<0}z^{-\rho/\ell}\Gamma\Bigl(\frac{\rho}{\ell}\Bigr) \Big \vert^2 \dx \alpha + N^{1/\ell}\xi
 \end{equation}
if $\xi>1/N$. We will treat only the first integral on the right hand
side of \eqref{R1-estim}, the second being completely similar. 
Clearly
\begin{equation}
\label{split-int}
\int_{1/N}^{\xi}
\Big \vert\sum_{\gamma>0}z^{-\rho/\ell}\Gamma\Bigl(\frac{\rho}{\ell}\Bigr) \Big \vert^2 
\dx \alpha 
= 
\sum_{k=1}^K
\int_\eta^{2\eta} \Big \vert\sum_{\gamma>0}z^{-\rho/\ell}\Gamma\Bigl(\frac{\rho}{\ell}\Bigr) \Big \vert^2 \dx \alpha
 +
\Odi{1} 
 \end{equation}
where $\eta=\eta_k= \xi/2^k$, $1/N\le \eta \le \xi/2$  and $K$ is a suitable integer satisfying $K=\Odi{L}$. 
Writing $\arctan 2\pi N\alpha = \pi/2 - \arctan(1/2\pi N\alpha)$ and
using the Saffari-Vaughan technique we have
\begin{align}
\notag
\int_{\eta}^{2\eta} 
\Big \vert\sum_{\gamma>0}z^{-\rho/\ell}\Gamma\Bigl(\frac{\rho}{\ell}\Bigr) 
\Big \vert^2\ \dx \alpha 
&\le 
\int_1^2 
\Bigl(
\int_{\delta\eta/2}^{2\delta\eta}
\Big \vert\sum_{\gamma>0}z^{-\rho/\ell}\Gamma\Bigl(\frac{\rho}{\ell}\Bigr) \Big \vert^2 \dx \alpha 
\Bigr) 
\dx \delta  
 \\
\label{expl}
&= \sum_{\gamma_1>0}\sum_{\gamma_2>0}
\Gamma\Bigl(\frac{\rho_{1}}{\ell}\Bigr)
\overline{\Gamma\Bigl(\frac{\rho_{2}}{\ell}\Bigr)}\ 
e^{\frac{\pi}{2\ell}(\gamma_1+\gamma_2)}\, \cdot  J,
\end{align}
say, where 
\[
J = J(N,\eta,\gamma_1,\gamma_2) 
= 
\int_1^2
\Bigl(
\int_{\delta\eta/2}^{2\delta\eta}
f_{1}(\alpha)f_2(\alpha)\ \dx \alpha 
\Bigr)
\dx \delta, 
\quad
w= \frac{1}{\ell} + \frac{i}{\ell}(\gamma_1-\gamma_2),
\]
\[
f_{1}(\alpha)=  
\vert z \vert^{-w} 
\quad \text{and}\quad 
f_{2}(\alpha) = 
\exp\Bigl(-\frac{\gamma_1+\gamma_2}{\ell}\arctan\frac{1}{2\pi N\alpha}\Bigr).
\]

Now we proceed to the estimation of $J$. Integrating twice by parts and
denoting by $F_1$ a primitive of $f_1$ and by $G_1$ a primitive of
$F_1$, we get
\begin{align}
\notag
J
&= 
\frac{1}{2\eta}
\Bigl(G_1(4\eta)f_2(4\eta)-G_1(2\eta)f_2(2\eta) \Bigr) 
-
\frac{2}{\eta}
\Bigl( G_1(\eta)f_2(\eta) - G_1\Bigl(\frac{\eta}{2}\Bigr)f_2\Bigl(\frac{\eta}{2}\Bigr) \Bigr) 
\\
\label{twiceparts}
&-
2\int_1^2  \! \! \!
G_1(2\delta\eta)f_{2}^{\prime}(2\delta\eta) \dx \delta 
+ 
2\int_1^2  \! \! \! \!
G_1\Bigl(\frac{\delta \eta}{2}\Bigr)f_{2}^{\prime}\Bigl(\frac{\delta \eta}{2}\Bigr) \dx \delta  
+
\int_1^2 \! \!
\Bigl(
\int_{\delta\eta/2}^{2\delta\eta}
G_1(\alpha)f_{2}^{\second}(\alpha)\, \dx \alpha 
\Bigr)
\dx \delta.
\end{align}
If $\alpha>1/N$ we have
\begin{align*}
f_{2}^{\prime}(\alpha) 
&
\ll_{\ell} \frac{1}{\alpha}
\Bigl(\frac{\gamma_1+\gamma_2}{N\alpha}\Bigr) f_{2}(\alpha)
\\
f_{2}^{\second}(\alpha) 
&\ll_{\ell} \frac{1}{\alpha^2}
\Bigl\{\Bigl(\frac{\gamma_1+\gamma_2}{N\alpha}\Bigr) +
\Bigl(\frac{\gamma_1+\gamma_2}{N\alpha}\Bigr)^2 \Bigr\} f_2(\alpha), 
\end{align*}
hence from \eqref{twiceparts} we get
\begin{equation} 
\label{J-estim}
J \ll_{\ell} \frac{1}{\eta}\max_{\alpha\in[\eta/2,4\eta]} \vert G_1(\alpha) \vert
\Bigl\{1+\Bigl(\frac{\gamma_1+\gamma_2}{N\eta}\Bigr)^2\Bigr\} 
\exp\Bigl(-c\Bigl(\frac{\gamma_1+\gamma_2}{N\eta}\Bigr)\Bigr),
\end{equation}
where $c=c(\ell)>0$ is a suitable constant.

In order to estimate $G_1(\alpha)$ we use the substitution 
\begin{equation}
\label{subst}
u = u(\alpha) = \Bigl(\frac{1}{N^2}+4\pi^2\alpha^2\Bigr)^{1/2},
\end{equation}
thus getting 
\[
F_1(\alpha) = 
\frac{1}{2\pi}
\int \frac{u^{1-w}}{(u^2 - N^{-2})^{1/2}}\,\dx u .
\]
By partial integration we have
\begin{equation}
\label{first-int}
F_1(\alpha) 
=
\frac{1}{2\pi(2-w)}
\Bigl\{ 
\frac{u^{2-w}}{(u^2 - N^{-2})^{1/2}} 
+
\int 
\frac{u^{3-w}}{(u^2 - N^{-2})^{3/2}}\, \dx u
\Bigr\}. 
\end{equation}
From \eqref{subst} and \eqref{first-int} we get
\begin{equation}
\label{twice-int}
G_1(\alpha) = 
\frac{1}{2\pi(2-w)}
 \Bigl\{
 A(\alpha) +
  \int B(\alpha)\, \dx \alpha
   \Bigr\},
\end{equation}
where 
\[
A(\alpha) = 
\frac{1}{2\pi} 
\int \frac{u^{3-w}}{u^2 - N^{-2}}\, \dx u
\quad
\textrm{and}
\quad
B(\alpha) = 
\int \frac{u^{3-w}}{(u^2 - N^{-2})^{3/2}}
\ \dx u.
\]
Again by partial integration we obtain 
\[
A(\alpha)  = 
\frac{1}{2\pi(4-w)}  
\Bigl\{
\frac{u^{4-w}}{u^2 - N^{-2}} 
+
2 \int \frac{u^{5-w}}{(u^2 - N^{-2})^2}\dx u
\Bigr\} 
\]
and 
\[
B(\alpha) = 
\frac{1}{4-w}  
\Bigl\{
\frac{u^{4-w}}{(u^2 - N^{-2})^{3/2}} 
+
3 \int \frac{u^{5-w}}{(u^2 - N^{-2})^{5/2}}\, \dx u
\Bigr\}. 
\]
Hence by \eqref{subst} we have for $\alpha \in [\eta/2,4 \eta]$ that 
\begin{align}
\label{A-B-estim}
A(\alpha) 
\ll_{\ell} 
\frac{u^{2-1/\ell}}{1+ \vert\gamma_1-\gamma_2 \vert} 
\ll
\frac{\alpha^{2-1/\ell}}{1+ \vert\gamma_1-\gamma_2 \vert} 
\quad
\textrm{and}
\quad
B(\alpha) 
\ll_{\ell}
\frac{\alpha^{1-1/\ell}}{1+ \vert\gamma_1-\gamma_2 \vert},
\end{align}
where $A(\alpha)$ and $B(\alpha)$ satisfy $A(\eta/4)=B(\eta/4)=0 $, and from
\eqref{twice-int}-\eqref{A-B-estim} we obtain  
\begin{equation}
\label{G1-estim}
G_1(\alpha) 
\ll_{\ell}  
\frac{\alpha^{2-1/\ell}}{1+ \vert\gamma_1-\gamma_2 \vert^2}
\end{equation}
for $\alpha\in[\eta/2,4\eta]$.
From \eqref{J-estim} and \eqref{G1-estim} we get
\[
J 
\ll_{\ell}
\eta^{1-1/\ell} 
\frac{1+\ (\frac{\gamma_1+\gamma_2}{N\eta})^2} {1+ \vert\gamma_1-\gamma_2 \vert^2}
\exp\Bigl(-c\Bigl(\frac{\gamma_1+\gamma_2}{N\eta}\Bigr)\Bigr),
\]
hence from \eqref{expl} and Stirling's formula we have
\begin{align}
\notag
\int_{\eta}^{2\eta} 
&\Big \vert
\sum_{\gamma>0}z^{-\rho/\ell}\Gamma\Bigl(\frac{\rho}{\ell}\Bigr) 
\Big \vert^2
\dx \alpha 
\\
\label{int-estim}
&\ll_{\ell} 
\eta^{1-1/\ell} 
\sum_{\gamma_1>0}\sum_{\gamma_2>0}
\vert \gamma_1 \vert^{(1-\ell)/(2\ell)} \vert \gamma_2 \vert^{(1-\ell)/(2\ell)}
\frac{1+(\frac{\gamma_1+\gamma_2}{N\eta})^2}{1+ \vert\gamma_1-\gamma_2 \vert^2}
\exp\Bigl(-c\Bigl(\frac{\gamma_1+\gamma_2}{N\eta}\Bigr)\Bigr) .
\end{align}
But sorting imaginary parts it is clear that
\[
\vert \gamma_1 \vert^{(1-\ell)/(2\ell)}
\vert \gamma_2 \vert^{(1-\ell)/(2\ell)}
\Bigl\{
1+\Bigl(\frac{\gamma_1+\gamma_2}{N\eta}\Bigr)^2
\Bigr\}
\exp\Bigl(-c\Bigl(\frac{\gamma_1+\gamma_2}{N\eta}\Bigr)\Bigr) 
\ll_{\ell}
 \vert \gamma_1 \vert^{(1-\ell)/\ell}
\exp\Bigl(-\frac{c}{2}\frac{\gamma_1}{N\eta}\Bigr),
\]
hence \eqref{int-estim} becomes
\begin{equation}
\label{int-estim1}
\ll_{\ell} \eta^{1-1/\ell} \sum_{\gamma_1>0} \vert \gamma_1 \vert^{(1-\ell)/\ell} \exp\Bigl(-\frac{c}{2}\frac{\gamma_1}{N\eta}\Bigr)
\sum_{\gamma_2>0}\frac{1}{1+ \vert\gamma_1-\gamma_2 \vert^2} \ll_{\ell} N^{1/\ell}\eta L^2,
\end{equation}
since the number of zeros $\rho_2= 1/2+i\gamma_2$ with $n\le  \vert\gamma_1-\gamma_2 \vert\le n+1$
is $\Odi{\log (n+ \vert\gamma_1 \vert)}$.

From \eqref{R1-first}-\eqref{split-int} and \eqref{int-estim1} we get
\begin{equation}
\label{last}
\int_{-\xi}^{\xi} 
\Big \vert
\sum_{\gamma>0}z^{-\rho/\ell}\Gamma\Bigl(\frac{\rho}{\ell}\Bigr) 
\Big \vert^2  
\dx \alpha 
\ll_{\ell} 
N^{1/\ell}\xi L^2,
\end{equation}
and  Lemma \ref{LP-Lemma-gen}  follows from
\eqref{last}.
\end{Proof}

We will also need the following result based on the Laplace formula 
for the Gamma function, see \cite{Laplace1812}.
In fact we will need it just for $\mu=2$ but, as before, we 
write the more general case.
\begin{Lemma}
\label{Laplace-formula}
Let $N$ be a positive integer, $z = 1/N - 2 \pi i \alpha$, and 
$\mu > 0$.
Then
\[
  \int_{-1 / 2}^{1 / 2} z^{-\mu} e(-n \alpha) \, \dx \alpha
  =
  e^{- n / N} \frac{n^{\mu - 1}}{\Gamma(\mu)}
  +
  \Odip{\mu}{\frac{1}{n}},
\]
uniformly for $n \ge 1$.
\end{Lemma}
\begin{Proof}
We start with the identity
\[
  \frac1{2 \pi}
  \int_{\R} \frac{e^{i D u}}{(a + i u)^s} \, \dx u
  =
  \dfrac{D^{s - 1} e^{- a D}}{\Gamma(s)},
\]
which is valid for $\sigma = \Re(s) > 0$ and $a \in \C$ with
$\Re(a) > 0$ and $D > 0$.
Letting $u = -2 \pi \alpha$ 
and taking $s = \mu$, $D = n$ and $a = N^{-1}$ we find
\[
  \int_{\R} \frac{e(- n \alpha)}{(N^{-1} - 2 \pi i \alpha)^\mu} \, \dx \alpha
  =
  \int_{\R} z^{- \mu} e(- n \alpha) \, \dx \alpha
  =
  \dfrac{n^{\mu - 1} e^{- n / N}}{\Gamma(\mu)}.
\]

For $0 < X < Y$ let
\[
  I(X, Y)
  =
  \int_X^Y \frac{e^{i D u}}{(a + i u)^\mu} \, \dx u.
\]
An integration by parts yields
\[
  I(X, Y)
  =
  \Bigl[
    \frac1{i D} \, \frac{e^{i D u}}{(a + i u)^\mu}
  \Bigr]_X^Y
  +
  \frac{\mu}D
  \int_X^Y \frac{e^{i D u}}{(a + i u)^{\mu + 1}} \, \dx u.
\]
Since $a > 0$, the first summand is $\ll_{\mu} D^{-1} X^{-\mu}$, uniformly.
The second summand is
\[
  \ll 
  \frac{\mu}D
  \int_X^Y \frac{\dx u}{u^{\mu + 1}}
  \ll_\mu
  D^{-1} X^{-\mu}.
\]
The result follows.
\end{Proof}

We remark that if $\mu\in \N$, $\mu\ge 2$,  Lemma~\ref{Laplace-formula} 
can be proved in an easier way  using the Residue Theorem 
(see, \emph{e.g.}, Languasco \cite{Languasco2000a} or 
Languasco and Zaccagnini  \cite{LanguascoZ2012b}).  

In the following we will also need a fourth-power average of $\Stilde_2(\alpha)$.

\begin{Lemma}
\label{Hua-Rieger-Lemma}
We have 
\[
\int_{-1/2}^{1/2} \,
\vert
\Stilde_2(\alpha) 
\vert^{4}
\, \dx \alpha 
\ll
N L^{2}.
\]
\end{Lemma} 
\begin{Proof}
Let $\mathcal{P}^{2}= \{p^{j}: j\ge 2, p\ \textrm{prime}\}$ and $r_{0}(m)$ be the number of representations of $m$ as a sum of two squares.
We have
\begin{align}
\notag
\int_{-1/2}^{1/2} & \,
\vert
\Stilde_2(\alpha) 
\vert^{4}\,
\dx \alpha 
\\
\notag
&=\!\!\!\!\!\!
\sum_{n_{1},n_{2},n_{3},n_{4}\ge 2}\!\!\!\!
\Lambda(n_{1})\Lambda(n_{2})\Lambda(n_{3})\Lambda(n_{4})\,
e^{-(n_{1}^{2}+n_{2}^{2}+n_{3}^{2}+n_{4}^{2})/N}
\int_{-1/2}^{1/2} \!\!\!\!
e((n_{1}^{2}+n_{2}^{2}-n_{3}^{2}-n_{4}^{2})\alpha)  
\, \dx \alpha 
\\
\notag
&
\ll
\sum_{p_{1},p_{2}\ge 2}
\log p_{1}\log p_{2}\,
e^{-2(p_{1}^{2}+p_{2}^{2})/N}
\sum_{\substack{p_{3},p_{4}\ge 2 \\ p_{1}^{2}+p_{2}^{2}=p_{3}^{2}+p_{4}^{2}}}
\log p_{3}\log p_{4}  
\\
\notag
&
\hskip1cm +
\sum_{\substack{n_{1},n_{2}\ge 2 \\n_{1}\in \mathcal{P}^{2}}}
\Lambda(n_{1})\Lambda(n_{2})\,
e^{-2(n_{1}^{2}+n_{2}^{2})/N}
\sum_{\substack{n_{3},n_{4} \ge 2 \\ n_{1}^{2}+n_{2}^{2}=n_{3}^{2}+n_{4}^{2}}}
\Lambda(n_{3})\Lambda(n_{4})
\\
\label{S-split}
&=\Sigma_{1}+\Sigma_{2},
\end{align}
say.
In $\Sigma_{1}$ we separately consider the contribution of the cases
$p_{1}p_{2} =p_{3} p_{4}$ and $p_{1}p_{2} \neq p_{3} p_{4}$; hence 
\(
\Sigma_{1}  \ll S_{1} + S_{2}
\)
where, by partial summation and  the Prime Number Theorem, we have
\[
S_{1}
=
2
\Bigl(
\sum_{p \ge 2}
(\log p )^{2} 
e^{-2p^{2}/N} 
\Bigr)^{2}
\ll
\Bigl(
1+
\int_2^{+\infty} \frac{u^2}{N}(\log u)\, e^{-2u^{2}/N}\, \dx u
\Bigr)^{2}
\ll 
N L^{2},
\]
and, by a dissection argument and  Satz 3 on page 94 of Rieger \cite{Rieger1968},
we also obtain
 \begin{align*}
S_{2}  &\ll 
\sum_{y\ge 1}\sum_{1\le x\le y}
y^{2} x^{2} 
e^{-2^{2y+1}/N}e^{-2^{2x+1}/N}
\Bigl(
\sum_{2^y\le p_{1}< 2^{y+1}} 
\sum_{2^x\le p_{2}< 2^{x+1}} 
\sum_{\substack{p_{3},p_{4} \ge 2\\ p_{1}^{2}+p_{2}^{2}=p_{3}^{2}+p_{4}^{2}\\ p_{1}p_{2} \neq p_{3} p_{4}}} 1
\Bigr)
\\
&
\ll
\sum_{y\ge 1} y^{4}   e^{-2^{2y+1}/N}
\Bigl(
\sum_{p_{1}, p_{2} < 2^{y+1}}  
\sum_{\substack{p_{3},p_{4} \ge 2\\ p_{1}^{2}+p_{2}^{2}=p_{3}^{2}+p_{4}^{2}\\ p_{1}p_{2} \neq p_{3} p_{4}}} 1
\Bigr)
\Bigl(
\sum_{1\le x\le y} e^{-2^{2x+1}/N}
\Bigr)
\\
&
\ll 
\sum_{y\ge 1}
y 2^{y} e^{-2^{2y+1}/N} 
\Bigl(
\int_{1}^{y} e^{-2^{t}/N} \dx t
\Bigr) 
\ll 
\sum_{y\ge 1}
y^{2} 2^{y} e^{-2^{2y+1}/N} 
\\
&
 \ll 
 \int_2^{+\infty} (\log u)^2 e^{-u/N}\, \dx u
\ll
N L^{2}.
\end{align*}
Summing up
\begin{equation}
\label{S1-estim}
\Sigma_{1}
 \ll  
N L^{2}.
\end{equation}
Recalling that $r_{0}(m) \ll m^{\eps}$,
it is also easy to see that
\begin{align}
\notag
\Sigma_{2}
&\ll
\sum_{\substack{n_{1},n_{2}\ge 2 \\n_{1}\in \mathcal{P}^{2}}}
\Lambda(n_{1})\Lambda(n_{2})(\log (n_{1}^{2}+n_{2}^{2}))^{2} \,
r_{0}(n_{1}^{2}+n_{2}^{2})\,
e^{-2(n_{1}^{2}+n_{2}^{2})/N}  
\\
\notag
&\ll
\sum_{\substack{n_{1},n_{2}\ge 2 \\n_{1}\in \mathcal{P}^{2}}}
n_{1}^{\eps}n_{2}^{\eps}\,
e^{-2(n_{1}^{2}+n_{2}^{2})/N}  
\ll
\Bigl(
\sum_{j\ge 2}
\sum_{p\ge 2} p^{j\eps} e^{-2p^{2j}/N}  
\Bigr)
\Bigl(
\sum_{n\ge 2} n^{\eps} e^{-2n^{2}/N}  
\Bigr)
\\
\notag
&
\ll
\Bigl(
\sum_{j\ge 2}  e^{-2^{2j}/N}  
\int_{2}^{+\infty} t^{j\eps } e^{-t^{2j}/N} \dx t  
\Bigr)
\Bigl(
N^{1/2+\eps} \int_{0}^{+\infty} u^{\eps-1/2} e^{-u}\, \dx u
\Bigr)
  \\ 
\notag
&
\ll
N^{1/2+2\eps} 
\sum_{j\ge 2}   N^{1/(2j)}  e^{-2^{2j}/N}   
  \\ 
\label{S2-estim} 
&
\ll
N^{1/2+2\eps} 
\Bigl(
N^{1/4} \log N 
+
\sum_{j > (1/2)\log N}     e^{-2^{2j}/N} 
\Bigr)
\ll
N^{3/4+3\eps}.
\end{align}

Combining \eqref{S-split}-\eqref{S2-estim}, Lemma \ref{Hua-Rieger-Lemma} follows.

\end{Proof}

\section{Proof of Theorem \ref{existence}}

Let $H\ge 2$, $H=\odi{N}$  be an  integer. 
We recall that we set $L = \log N$ for brevity.
Recalling \eqref{r-def} and letting 
\[
  R(n)
  =
  \sum_{a + b^2 + c^{2} = n} \Lambda(a)\Lambda(b)\Lambda(c),
\]
we have (see, \emph{e.g.}, page 14 of \cite{Zhao2014a}) that 
\begin{equation}
\label{R-r-diff}
r(n) = R(n) + \Odi{n^{3/4} (\log n)^{3}}.
\end{equation}
Then, for every $n\le 2N$, we can write
\[
  r(n)
  =
   R(n) +  \Odi{n^{3/4}  (\log n)^{3}}
  =
  e^{n / N} 
  \int_{-1/2}^{1/2}
    \Stilde_1(\alpha)  \Stilde_2(\alpha)^{2}  e(-n \alpha) \, \dx \alpha
   + \Odi{n^{3/4}  (\log n)^{3}}.
\]
From this equation, the Cauchy-Schwarz inequality, Lemma \ref{Hua-Rieger-Lemma} and the Prime Number Theorem,
for every $n\le 2N$
we also have
\begin{align}
 \notag
  r(n)
 &\ll 
  \int_{-1/2}^{1/2}
    \vert \Stilde_1(\alpha) \vert \vert  \Stilde_2(\alpha) \vert^{2}  \dx \alpha 
    + N^{3/4}L^{3}
    \\
    \label{upperbound}
&   \ll
    \Bigl(
     \int_{-1/2}^{1/2}
    \vert \Stilde_1(\alpha) \vert ^{2}  \dx \alpha
    \Bigr)^{1/2}
    \Bigl(
    \int_{-1/2}^{1/2}
    \vert  \Stilde_2(\alpha) \vert^{4}  \dx \alpha
    \Bigr)^{1/2}
    +N^{3/4}L^{3}
    \ll
    N L^{3/2}.
\end{align}

We need now to choose a suitable weighted average of $r(n)$.
We further set
\begin{equation}
\notag 
   U(\alpha,H)
  = 
  \sum_{m=1}^H e(m \alpha)  
\end{equation}
and, moreover, we also have the usual 
numerically explicit inequality
\begin{equation}
\label{UH-estim}
\vert U(\alpha,H) \vert
\le
\min \Bigl(H; \frac{1}{\vert \alpha\vert} \Bigr).
\end{equation}
With these definitions and \eqref{R-r-diff}, we may write
\begin{align*}
  \Stilde (N, H)
  &: =
  \sum_{n = N+1}^{N + H} 
  e^{-n / N} r(n)
  =
  \int_{-1/2}^{1/2}
    \Stilde_1(\alpha) \Stilde_2(\alpha)^{2}  U(-\alpha,H) e(-N\alpha) \, \dx \alpha
   +
   \Odi{H N^{3/4}L^{3}}.
\end{align*} 

Using Lemma \ref{Linnik-lemma} with $\ell=1,2$ 
and recalling that
$\Gamma(1)=1$, $\Gamma(1/2)=\pi^{1/2}$, we can write
\begin{align}
\notag
  \Stilde (N, H)
  &= 
    \int_{-1/2}^{1/2}
    \frac{\pi}{4z^{2}} U(-\alpha,H) e(-N\alpha) \, \dx \alpha
    +
       \int_{-1/2}^{1/2}
       \frac{1}{z} \Bigl( \Stilde_2(\alpha)^{2}-\frac{\pi}{4z}\Bigr) 
    U(-\alpha,H) e(-N\alpha) \, \dx \alpha
\\
\notag
&+
  \int_{-1/2}^{1/2}
    \Bigl( \Stilde_1(\alpha) -\frac{1}{z}\Bigr)\Stilde_2(\alpha)^{2} U(-\alpha,H) e(-N\alpha) \, \dx \alpha 
+
    \Odi{H N^{3/4}L^{3}}
    \\
    \label{main-split}
    & =
    I_{1}+I_{2}+I_{3} +  \Odi{H N^{3/4}L^{3}},
\end{align}
say. From now on, we denote 
\[
\Etilde_{\ell}(\alpha) : =\Stilde_\ell(\alpha) - \frac{\Gamma(1/\ell)}{\ell z^{1/\ell}}.
\]

Using Lemma~\ref{Laplace-formula} we immediately get 
\begin{equation}
\label{I1-eval}
I_{1}
=
\frac{\pi}{4}   \sum_{n = N+1}^{N + H}  ne^{-n/N} + \Odi{\frac{H}{N}}
  =
\frac{\pi H  N}{4e} + \Odi{H^{2}}.
\end{equation}

Now we estimate $I_{2}$.
Using the identity $f^{2}-g^{2} = 2f(f-g) -(f-g)^{2}$ we obtain 
\begin{equation}
\label{I2-estim1}
I_{2} 
\ll
 \int_{-1/2}^{1/2}
      \vert \Etilde_{2}(\alpha)  \vert
        \frac{\vert U(\alpha,H)\vert}{\vert z\vert^{3/2}}
     \, \dx \alpha
+
 \int_{-1/2}^{1/2}
       \vert \Etilde_{2}(\alpha)  \vert^{2}  
       \frac{\vert U(\alpha,H) \vert}{\vert z\vert} 
        \, \dx \alpha 
    =
   J_{1}+J_{2},
\end{equation}
say.
Using   \eqref{z-estim},  \eqref{UH-estim}, Lemma \ref{LP-Lemma-gen} and a partial integration argument we obtain
\begin{align}
\notag
J_{2}
&\ll 
HN
 \int_{-1/N}^{1/N}  \vert \Etilde_{2}(\alpha)  \vert^{2}   \, \dx \alpha 
 +
 H \int_{1/N}^{1/H}  \vert \Etilde_{2}(\alpha)  \vert^{2}  
 \,\frac{\dx \alpha}{\alpha} 
 +
 \int_{1/H}^{1/2}  \vert \Etilde_{2}(\alpha)  \vert^{2}   
 \,\frac{\dx \alpha}{\alpha^{2}}
 \\
 \notag
 &\ll 
HN^{1/2} L^{2} 
+ 
H N^{1/2} L^{2}\Bigl(1  + \int_{1/N}^{1/H} \frac{\dx \xi}{\xi}   \Bigr)
+
N^{1/2} L^{2} \Bigl(H +   \int_{1/H}^{1/2}   \frac{\dx \xi}{\xi^2} \Bigr)
 \\
 \label{J2-estim}
 &\ll 
H N^{1/2} L^{3}. 
\end{align} 
Using  the Cauchy-Schwarz inequality and arguing as for $J_{2}$ we get 
\begin{align}
\notag
J_{1}
&\ll 
H N^{3/2} \Bigl(  \int_{-1/N}^{1/N}\!\!\!\! \dx \alpha  \Bigr)^{1/2}
\Bigl(
 \int_{-1/N}^{1/N} \!\! \vert \Etilde_{2}(\alpha)  \vert^{2}   \, \dx \alpha 
 \Bigr)^{1/2}
 +
 H \Bigl(  \int_{1/N}^{1/H} \!\! \frac{\dx \alpha}{\alpha^{2}}\Bigr)^{1/2} 
\Bigl(
 \int_{1/N}^{1/H}\!\!  \vert \Etilde_{2}(\alpha)  \vert^{2}  
 \frac{\dx \alpha}{\alpha}
  \Bigr)^{1/2}
  \\
  \notag
     &\hskip2cm 
     +
  \Bigl(  \int_{1/H}^{1/2}   \frac{\dx \alpha}{\alpha^{4}}\Bigr)^{1/2}  
\Bigl(
 \int_{1/H}^{1/2}  \vert \Etilde_{2}(\alpha)  \vert^{2}   
 \frac{\dx \alpha}{\alpha}
   \Bigr)^{1/2}
 \\
 \notag
 &\ll 
H N^{3/4} L 
+ 
H N^{3/4} L  \Bigl( 1 + \int_{1/N}^{1/H}  \frac{\dx \xi}{\xi}   \Bigr)^{1/2} 
 +
H^{3/2} N^{1/4} L \Bigl( 1+   \int_{1/H}^{1/2}  \frac{\dx \xi}{\xi}   \Bigr)^{1/2}
 \\
 \label{J1-estim}
 &\ll 
H N^{3/4} L^{3/2} + H^{3/2}N^{1/4} L^{3/2}
\ll
H N^{3/4} L^{3/2}. 
\end{align}
Combining \eqref{I2-estim1}-\eqref{J1-estim} we finally obtain
\begin{equation}
\label{I2-estim-final}
I_{2} \ll H N^{3/4} L^{3/2} .
\end{equation}

Now we estimate $I_{3}$.
By the Cauchy-Schwarz  inequality,  \eqref{UH-estim} and Lemma
\ref{Hua-Rieger-Lemma}  we obtain
\begin{align}
\notag
I_{3}
&\ll 
\Bigl( \int_{-1/2}^{1/2}  \vert\Stilde_2(\alpha)\vert^{4}   \, \dx \alpha \Bigr)^{1/2}
\Bigl( 
\int_{-1/2}^{1/2}  \vert \Etilde_{1}(\alpha)  \vert^{2} \vert U(\alpha,H)\vert^2 
 \, \dx \alpha 
 \Bigr)^{1/2}
\\
\notag
&\ll
N^{1/2}L
\Bigl( H^2
\int_{-1/H}^{1/H}  \vert \Etilde_{1}(\alpha)  \vert^{2} \, \dx \alpha
+ 
 \int_{1/H}^{1/2}  \vert \Etilde_{1}(\alpha)  \vert^{2}   \, \frac{\dx \alpha}{\alpha^{2}} 
 \Bigr)^{1/2}
 \\
 \label{I3-estim-final}
 &\ll
H^{1/2}NL^2,
 \end{align}
where in the last step we used  Lemma \ref{LP-Lemma-gen} and 
a partial integration argument.

By \eqref{main-split}-\eqref{I1-eval}, \eqref{I2-estim-final} and \eqref{I3-estim-final},
we can finally write
\begin{equation*}
\Stilde (N, H) = \frac{\pi}{4e} H N
+   
\Odi{H^{1/2}N L^{2}+H N^{3/4} L^{3} + H^{2}}.
\end{equation*}
Theorem \ref{existence} follows since the exponential weight $e^{-n/N}=e^{-1}+ \Odi{H/N}$ for $n\in[N+1,N+H]$ and hence by \eqref{upperbound} 
it can be removed at the cost of 
inserting an extra factor $\Odi{H^{2}L^{3/2}}$ in the error term. 
 The corollary about the existence in short intervals follows by remarking that  $\Stilde (N, H) > 0$ if $L^{4}\ll H  = \odi{NL^{-3/2}}$.
\qed

\vskip0.5cm
\noindent
\begin{tabular}{l@{\hskip 18mm}l}
Alessandro Languasco               & Alessandro Zaccagnini\\
Universit\`a di Padova     & Universit\`a di Parma\\
Dipartimento di Matematica & Dipartimento di Matematica e Informatica \\
Via Trieste 63                & Parco Area delle Scienze, 53/a \\
35121 Padova, Italy            & 43124 Parma, Italy\\
{\it e-mail}: languasco@math.unipd.it        & {\it e-mail}:
alessandro.zaccagnini@unipr.it  
\end{tabular}


\begin{thebibliography}{10}
\providecommand{\url}[1]{\texttt{#1}}
\providecommand{\urlprefix}{URL }
\expandafter\ifx\csname urlstyle\endcsname\relax
  \providecommand{\doi}[1]{doi:\discretionary{}{}{}#1}\else
  \providecommand{\doi}{doi:\discretionary{}{}{}\begingroup
  \urlstyle{rm}\Url}\fi

\bibitem{HardyL1923}
G.~H. Hardy, J.~E. Littlewood -  \textsl{Some problems of `{P}artitio
  numerorum'; {III}: {O}n the expression of a number as a sum of primes} -
  Acta Math., \textbf{44} (1923), 1--70.
  
  \bibitem{HarmanK2010}
G.~Harman, A.~Kumchev -  \textsl{{On sums of squares of primes II}} -  J. Number Theory, \textbf{130} (2010), 1969--2002.

\bibitem{Hua1938}
L.~K. Hua -  \textsl{Some results in the additive prime number theory} -
  Quart. J. Math. Oxford, \textbf{9} (1938), 68--80.

\bibitem{Languasco2000a}
A.~Languasco -  \textsl{Some refinements of error terms estimates for certain
  additive problems with primes} -  J. Number Theory, \textbf{81} (2000),
  149--161.

\bibitem{LanguascoP1994}
A.~Languasco, A.~Perelli -  \textsl{On {Linnik}'s theorem on {Goldbach} numbers
  in short intervals and related problems} -  Ann. Inst. Fourier, \textbf{44}
  (1994), 307--322.

\bibitem{LanguascoZ2012b}
A.~Languasco, A.~Zaccagnini -  \textsl{Sums of many primes} -  J. Number
  Theory, \textbf{132} (2012), 1265--1283.
  
\bibitem{LanguascoZ2013a}
A.~Languasco, A.~Zaccagnini -  \textsl{{A Ces\`aro Average of Hardy-Littlewood
  numbers}} -  J. Math. Anal. Appl., \textbf{401} (2013), 568--577.

\bibitem{LanguascoZ2013b}
A.~Languasco, A.~Zaccagnini -  \textsl{On a ternary Diophantine problem with
  mixed powers of primes} -  Acta Arith., \textbf{159} (2013), 345--362.

\bibitem{LanguascoZ2014a}
A.~Languasco, A.~Zaccagnini -  \textsl{{A Ces\`aro Average of Goldbach
  numbers}} -   Forum Math. \textbf{27} (2015), 1945--1960.

\bibitem{Laplace1812}
P.~S. Laplace -  \textsl{Th\'eorie analytique des probabilit\'es} -  Courcier
  (1812).

\bibitem{LeungL1993}
M.~Leung, M.~Liu -  \textsl{On generalized quadratic equations in three prime
  variables} -  Monatsh. Math., \textbf{115} (1993), 133--167.

\bibitem{Li2008}
H.~Li -  \textsl{{Sums of one prime and two prime squares}} -  Acta Arith., \textbf{134} (2008), 1--9.
  
\bibitem{Linnik1946}
Y.~V. Linnik -  \textsl{A new proof of the {Goldbach-Vinogradow} theorem} -
  Rec. Math. N.S., \textbf{19 (61)} (1946), 3--8, (Russian).

\bibitem{Linnik1952}
Y.~V. Linnik -  \textsl{Some conditional theorems concerning the binary
  {Goldbach} problem} -  Izv. Akad. Nauk SSSR Ser. Mat., \textbf{16} (1952),
  503--520, (Russian).
  

\bibitem{Pintz2012}
J.~Pintz -  \textsl{{The Bounded Gap Conjecture and bounds between consecutive
  Goldbach numbers}} -  Acta Arith., \textbf{155} (2012), 397--405.

\bibitem{Rieger1968}
G.~J. Rieger -  \textsl{\"{U}ber die {S}umme aus einem {Q}uadrat und einem
  {P}rimzahlquadrat} -  J. Reine Angew. Math., \textbf{231} (1968), 89--100.
  
\bibitem{Schwarz1961}
W.~Schwarz -  \textsl{Zur {D}arstellung von {Z}ahlen durch {S}ummen von
  {P}rimzahlpotenzen. {I}. {D}arstellung hinreichend grosser {Z}ahlen} -  J.
  Reine Angew. Math., \textbf{205} (1960/1961), 21--47.


\bibitem{Wang2004}
M.~Wang -  \textsl{On the sum of a prime and two prime squares} -  Acta Math.
  Sinica (Chin. Ser.), \textbf{47} (2004), 845--858.

\bibitem{WangM2006a}
M.~Wang, X.~Meng -  \textsl{{The exceptional set in the two prime squares and a
  prime problem}} -  Acta Math. Sinica (Eng. Ser.), \textbf{22}
  (2006), 1329--1342.
  
\bibitem{Zhao2014a}
L.~Zhao -  \textsl{{The additive problem with one prime and two squares of
  primes}} -  J. Number Theory, \textbf{135} (2014), 8--27.

\end{thebibliography}
\end{document}